\documentclass{conm-p-l}
\newtheorem{theorem}{Theorem}[section]
\newtheorem{lemma}[theorem]{Lemma}
\newtheorem{corollary}[theorem]{Corollary}
 
\theoremstyle{definition}  
\newtheorem{definition}[theorem]{Definition}
\newtheorem{example}[theorem]{Example}

\theoremstyle{remark}

\newcommand{\id}{\mbox{id}}

\newcommand{\End}{\mbox{End}} 

\newcommand{\eps}{\varepsilon} 

\newcommand{\act}{\rightharpoonup}

\newcommand{\la}{{<}}  
\newcommand{\ra}{{>}}  

\newcommand{\1}{_{(1)}} 
\newcommand{\2}{_{(2)}} 
\newcommand{\3}{_{(3)}} 
\newcommand{\4}{_{(4)}}

\begin{document}

\title{\bf A Duality Theorem for Quantum Groupoids}  

\author{Dmitri Nikshych} 
\address{Department of Mathematics, UCLA,  
405 Hilgard Avenue, Los Angeles, CA 90095-1555}
\email{nikshych@math.ucla.edu}

\date{}

\begin{abstract} 
We prove a duality theorem for quantum groupoid (weak Hopf algebra) 
actions that extends the well-known result for usual
Hopf algebras obtained in \cite{BM} and \cite{vdB}.  
\end{abstract}

\maketitle  

\begin{section}
{Introduction}

By {\em (finite) quantum groupoids} we understand weak Hopf algebras
introduced in \cite{BNSz}, \cite{BSz} as a generalization of ordinary
Hopf algebras providing a good framework for studying symmetries 
of certain  quantum field theories. These objects also generalize both 
ordinary  groupoid algebras  and their duals. 
A special case of quantum groupoids  with involutive antipode 
was studied in \cite{NV1}, \cite{N}. 

Finite quantum groupoids naturally arise in the theory of von Neumann
algebras : it was shown in \cite{NV2} that finite index II${}_{1}$ 
subfactors of depth $\leq 2$ can be characterized 
as $C^*$-quantum groupoid smash products. 
This result was extended in \cite{NV3}, where a uniform description
of all finite depth subfactors was obtained via a Galois correspondence.
In fact, one can use subfactors in order to construct interesting concrete 
examples of quantum groupoids such as Temperley-Lieb algebras 
\cite{NV2}, \cite{NV3}.

Another motivation to study quantum groupoids comes from the fact that
their representation theory provides examples of monoidal categories that 
can be used for constructing invariants of links and 3-manifolds \cite{NVT}.

In this paper we prove the following duality theorem for 
smashed products :  if $H$ is a finite quantum groupoid and $A$ is an
$H$-module algebra, then $(A\#H)\#H^* \cong \End(A\#H)_A$, where $H^*$
acts on $A\#H$ in a dual way and $A\#H$ is viewed as a right $A$-module
via multiplication. 
For usual Hopf algebras this result was proved in \cite{BM}
(where infinite dimensional case was considered) and \cite{vdB}. 
For weak Kac algebras (i.e., finite $C^*$-quantum groupoids 
with an involutive antipode) it was established in \cite{N}. 

The note is organized as follows.

In Preliminaries (Section $2$) we recall definitions and basic facts 
concerning finite quantum groupoids (weak Hopf algebras) and prove
identities we need for later computations.

In Section $3$ we prove  the main result by writing down explicit formulas
for isomorphism between $(A\#H)\#H^*$ and $\End(A\#H)_A$. As a corollary
we obtain that $H\#H^*$ is always a semisimple algebra.

The results of this paper were presented by the author at the
Colloquium on Quantum Groups and Hopf Algebras held in La Falda,
Argentina  in August, 1999 and he would like to thank N.~Andruskiewitsch,
W.~Ferrer Santos, and H.-J.~Schneider for inviting him. 
The author is also grateful to L.~Vainerman for numerous discussions
on quantum groupoids and his comments on the present work.

\end{section}


\begin{section} 
{Preliminaries}

Let $k$ be a field.

Throughout this paper we use Sweedler's notation for
comultiplication, writing $\Delta(b) = b\1 \otimes b\2$.

\begin{definition}
\label{basic definition}
By a {\em weak Hopf algebra} \cite{BNSz}, or {\em finite quantum groupoid}
we understand a finite dimensional $k$-vector space $H$
that has structures of algebra $(H,\,m,\,1)$ and coalgebra
$(H,\,\Delta,\,\eps)$ related as follows:
\begin{enumerate}
\item[(1)] $\Delta$ is a (not necessarily unit-preserving) homomorphism :
$$
\Delta(hg) = \Delta(h)\Delta(g), 
$$
\item[(2)] The unit and counit satisfy the identities
\begin{eqnarray*}
\eps(hgf) &=& \eps(hg\1)\eps(g\2f) =  \eps(hg\2)\eps(g\1f), \\
(\Delta \otimes \id) \Delta(1) &=& 
(\Delta(1)\otimes 1)(1\otimes \Delta(1)) =
(1\otimes \Delta(1))(\Delta(1)\otimes 1), 
\end{eqnarray*}
\item[(3)]
There is a linear map $S: H \to H$, called an {\em antipode}, such that
\begin{eqnarray*}
m(\id \otimes S)\Delta(h) &=&(\eps\otimes\id)(\Delta(1)(h\otimes 1)),\\
m(S\otimes \id)\Delta(h) &=& (\id \otimes \eps)((1\otimes h)\Delta(1)),\\
S(h\1)h\2 S(h\2) &=& S(h),
\end{eqnarray*}
\end{enumerate}
for all $h,g,f\in H$.
\end{definition}
The antipode is unique and invertible \cite{BNSz}, moreover
it is an anti-algebra and anti-coalgebra map. The right-hand
sides of two first formulas in $(3)$ are called {\em target}
and {\em source counital maps} and denoted $\eps_t$, $\eps_s$
respectively :
\begin{eqnarray*}
\eps_t(h) = (\eps\otimes\id)(\Delta(1)(h\otimes 1)),\\
\eps_s(h) = (\id \otimes \eps)((1\otimes h)\Delta(1)).
\end{eqnarray*}
The counital maps $\eps_t$ and $\eps_s$ are idempotents in $\End_k(H)$,
we also have relations $S\circ \eps_t = \eps_s \circ S$ and
$S\circ \eps_s = \eps_t \circ S$.

The main difference between quantum groupoids and Hopf algebras
is that the ranges of counital maps are, in general, separable
subalgebras of $H$ not necessarily equal to $k$. They are
called {\em target} and {\em source counital subalgebras} 
and play a role of ``non-commutative bases'' 
(cf. Example~\ref{examples} below) :
\begin{eqnarray*}
H_t &=& \{h\in H \mid \eps_t(h) =h \}
  =  \{h\in H \mid \Delta(h) = 1\1 h\otimes 1\2 = h 1\1 \otimes 1\2 \}, \\
H_s &=& \{h\in H \mid \eps_s(h) =h \} 
= \{h\in H \mid \Delta(h) = 1\1 \otimes h 1\2 =  1\1 \otimes 1\2 h \}.
\end{eqnarray*}  
The counital subalgebras commute, the restriction of the antipode
gives an anti-isomorphism between $B_t$ and $B_s$, moreover,
$B_t$ (resp.\ $B_s$) is a left (resp.\ right) coideal subalgebra 
of $B$. We also have $S\circ \eps_t = \eps_s \circ S$,
$S^2\vert_{B_t} =\id_{B_t}$, and $S^2\vert_{B_s} =\id_{B_s}$.

Note that $H$ is an ordinary Hopf algebra if and only if
$\Delta(1)=1\otimes 1$ if and only if $\eps$ is a homomorphism 
if and only if $H_t=H_s =k$.

The dual vector space $H^*$ has a natural structure of a quantum
groupoid with the structure operations dual to those of $H$ :
\begin{eqnarray*}
& & \la \phi\psi,\,h\ra = \la \phi\otimes\psi,\,\Delta(h) \ra, \\
& & \la \Delta(\phi),\,h \otimes g \ra =  \la \phi,\,hg \ra, \\
& & \la S(\phi),\,h\ra = \la \phi,\,S(h) \ra, 
\end{eqnarray*}
for all $\phi,\psi \in H^*,\, h,g\in H$. The unit of $H^*$ is $\eps$
and counit is $\phi \mapsto \la\phi,\, 1\ra$.

\begin{example}
\label{examples}
Let $G$ be a finite {\em groupoid} (a category with finitely many
morphisms, such that each morphism is invertible) then the groupoid
algebra $kG$ (generated by morphisms $g\in G$ with the product of 
two morphisms being equal to  their composition
if the latter is defined and $0$ otherwise) 
is a quantum groupoid via :
$$
\Delta(g) = g\otimes g,\quad \eps(g) =1,\quad S(g)=g^{-1},\quad g\in G.
$$ 
The dual quantum groupoid $(kG)^*$ is generated by idempotents
$p_g,\, g\in G$ such that $p_g p_h= \delta_{g,h}p_g$ and
$$
\Delta(p_g) =\sum_{uv=g}\,p_v\otimes p_v,\quad \eps(p_g)= \delta_{g,gg^{-1}},
\quad S(p_g) =p_{g^{-1}}.
$$
\end{example}
 
It is known that any group action on a set gives rise to a finite
groupoid \cite{R}. Similarly, in the ``quantum'' situation, one
can associate a weak Hopf algebra (quantum groupoid) with every
action of a usual Hopf algebra on a separable algebra, see \cite{NVT}
for details.

Finally, the most non-trivial examples of quantum groupoids known
so far come from the theory of von Neumann II${}_1$ subfactors \cite{GHJ} : 
in \cite{NV2}  finite index  subfactors of depth $\leq 2$ were characterized 
as quantum groupoid smash products and it was explained in \cite{NV3}  
that it is possible to construct concrete examples of quantum groupoids 
from subfactors of arbitrary finite depth.

An algebra $A$ is a left  {\em $H$-module algebra} \cite{NSzW}
if $A$ is a left $H$-module via $h\otimes x \mapsto h\cdot x$ and
$$
h\cdot (xy) = (h\1\cdot x)(h\2\cdot y), \qquad  h\cdot 1 =\eps_t(h)\cdot 1,
$$
for all $h\in H,\,x,y\in A$. 

A {\em smash product} algebra $A\#H$ of $A$ and $H$ is defined on a 
$k$-vector space $A \otimes_{H_t} H$ (relative tensor product), 
where $H$ is a left  $H_t$-module via multiplication and $A$ is a right 
$H_t$-module via 
$$
x\cdot z = S(z)\cdot x = x(z\cdot 1) \qquad x\in A, z\in H_t.
$$ 
Let $x\#h$ be the class of  $x\otimes h$ in  $A \otimes_{H_t} H$, 
then the multiplication of $A\#H$ is given by the familiar formula :
$$
(x\#h)(y\#g) = x(h\1 \cdot y)\#h\2 g, \qquad x,y\in A,\, h,g\in H
$$ 
and the unit of $A\#H$ is $1\#1$.
%

\begin{example}
\label{actions}
The target counital subalgebra $H_t$ is a {\em trivial $H$-module}
algebra with the action of $H$ given by $h\cdot z = \eps_t(hz)$,
where $h\in H,\, z\in H_t$.

The dual quantum groupoid $H^*$ is an $H$-module algebra via 
$$
h\act \phi = \phi\1 \la \phi\2,\, h\ra,
$$ 
for all $h\in H,\, \phi\in H^*$.
\end{example}

In the following Lemma we collect the identities we will
use in what follows. They can be found in \cite{BNSz} and 
\cite{NV1}, we include them here for the convenience 
of the reader.

\begin{lemma}
\label{a lemma}
For every quantum groupoid $H$ and elements $h\in H,\,z\in H_t$
the following identities hold true :
\begin{enumerate}
\item[(i)] $h\1 \otimes \eps_t(h\2)  =1\1h \otimes 1\2\quad$
and $\quad\eps_s(h\1) \otimes h\2 = 1\1 \otimes h 1\2$, 
\item[(ii)] $1\1 S(z) \otimes 1\2 = 1\1 \otimes 1\2z$,
\item[(iii)] $h\2 S^{-1}(h\1)\otimes h\3 
= S(\eps_t(h\1)\otimes h\2 = 1\1\otimes 1\2h$.
\end{enumerate}
\end{lemma}
\begin{proof}
(i) We have :
\begin{eqnarray*}
h\1 \otimes \eps_t(h\2)
&=& h\1 \eps(1\1h\2) \otimes 1\2 \\
&=& 1\1h\1 \eps(1\2h\2) \otimes 1\3 = 1\1h \otimes 1\2,
\end{eqnarray*}
where we used the definition of $\eps_t$ and the axiom (2)
of Definition~\ref{basic definition}. The second identity
is similar.
\newline(ii)
Since $S(z) \in H_s$ we can compute : 
\begin{eqnarray*}
1\1 S(z) \otimes 1\2 
&=& S(z)\1 \otimes \eps_t(S(z)\2) \\
&=& 1\1 \otimes \eps_t(1\2S(z)) \\
&=& 1\1 \otimes \eps_t(1\2z) = 1\1 \otimes 1\2z,
\end{eqnarray*}
using part (i), definition of the source counital subalgebra,
and the identity $\eps_t(hg) = \eps_t(h \eps_t(g))$
that follows from axiom (2) of Definition~\ref{basic definition}.
Observe that $S(1\1) \otimes 1\2$ is a  
separability idempotent \cite{P} of $H_t$.
\newline (iii)
Using part (ii) and the fact that $H_s$ and $H_t$ commute, we have
\begin{eqnarray*}
h\2 S^{-1}(h\1)\otimes h\3
&=& S(\eps_t(h\1)) \otimes h\2 \\
&=& S(\eps_t(1\1h\1)) \otimes 1\2h\2 \\
&=& 1\1 S(\eps_t(h\1)) \otimes 1\2h\2 \\
&=& 1\1 \otimes 1\2 \eps_t(h\1) h\2 = 1\1 \otimes 1\2h.
\end{eqnarray*}
\end{proof}
\end{section}


\begin{section}
{Main result}

Let $H$ be a finite quantum groupoid and
$A$ be a left $H$-module algebra. Then the smash product $A\#H$
is a left $H^*$-module algebra via
$$
\phi\cdot (a\#h) = a\#(\phi\act h), \quad
\phi\in H^*,\,h\in H,\,a\in A.
$$
In the case when $H$ is an ordinary finite dimensional Hopf algebra,
it follows from \cite{BM} that there is an isomorphism
$(A\#H)\#H^* \cong M_n(A)$, where $n=\dim H$ and $M_n(A)$ is an algebra
of $n$-by-$n$ matrices over $A$.

We will show that this result extends to quantum groupoid
action in the form $(A\#H)\#H^* \cong \End(A\#H)_A$, 
where $A\#H$ is a right $A$-module
via multiplication (note that $A\#H$ is not necessarily a free $A$-module,
so that we have $\End(A\#H)_A\not\cong M_n(A)$ in general;
see (\cite{NV2}, 7) for an example when $H$ is not free over $H_t$). 
We will explicitly  write down canonical isomorphisms 
between $(A\#H)\#H^*$ and $\End(A\#H)_A$.

\begin{lemma}
\label{alpha}
The map $\alpha : (A\#H)\#H^* \to \End(A\#H)_A$ defined by
$$
\alpha((x\#h)\#\phi)(y\#g) = (x\#h)(y\#(\phi\act g)) 
= x(h\1\cdot y) \#h\2 (\phi\act g) 
$$ 
for all $x,y\in A,\, h,g\in H,\,\phi\in H^*$ is a homomorphism
of algebras.
\end{lemma}
\begin{proof}
First, we need to check that $\alpha$ is well defined.
For all $z\in H_t$ and $\xi\in H_t^*$ we have :
\begin{eqnarray*}
\alpha((x\#zh)\#\phi)(y\#g)
&=& x(zh\1\cdot b)\#h\2(\phi\act g) \\
&=& x(z\cdot 1)(h\1\cdot b)\#h\2(\phi\act g) \\
&=& \alpha( (x\cdot z)\#h)\#\phi )(y\#g), \\
\alpha((x\#h)\#\xi\phi)(y\#g) 
&=& x(h\1\cdot b)\#h\2(\xi\act 1)(\phi\act g) \\
&=& \alpha( (x\#h(\xi\act 1))\#\phi )(y\#g) \\
&=&  \alpha( (x\#(S(\xi)\act h))\#\phi )(y\#g)\\
&=&  \alpha( (x\# (h\cdot \xi)) \# \phi )(y\#g)
\end{eqnarray*}
where we used definition of the target counital subalgebra,
Lemma~\ref{a lemma}(ii), and that $(\xi\act 1)\in H_s$ 
for all $\xi\in H_t^*$.

Next, we verify that $\alpha((x\#h)\#\phi) \in \End(A\#H)_A$
for all $x\in A,\,h\in H,\,\phi\in H^*$. For all $z\in H_t$
we have :
\begin{eqnarray*}  
\alpha((x\#h)\#\phi) (y\#zg)
&=& x(h\1\cdot y)\# h\2z(\phi\act g) \\
&=& x(h\1S(z) \cdot y)\# h\2 (\phi\act g) \\
&=& \alpha((x\#h)\#\phi)((y \cdot z) \# g),
\end{eqnarray*}
using the identity $\phi\act zg = z(\phi\act g)$ and 
Lemma~\ref{a lemma} (ii).

The following computation shows that $\alpha$
commutes with the right action of all $w\in A$ :
\begin{eqnarray*}
\alpha((x\#h)\#\phi) ((y\#g)\cdot w)
&=& \alpha((x\#h)\#\phi)(y(g\1\cdot w)\# g\2) \\
&=& (x\#h)(y(g\1\cdot w) \# (\phi\act g\2) \\
&=& (x\#h)(y\# (\phi\act g))(w\# 1) \\
&=& (\alpha((x\#h)\#\phi)(y\#g) )\cdot w.
\end{eqnarray*}
Finally,
\begin{eqnarray*}
\lefteqn{
\alpha(\, ((x\#h)\#\phi) ((x'\#h')\#\phi') \,)(y\#g) = } \\
&=& \alpha((x\#h)(x'\# (\phi\1\act h')) \# \phi\2\phi') (y\#g) \\
&=& (x\#h) (x'\# (\phi\1\act h')) (y\# (\phi\2\phi'\act g)) \\
&=& (x\#h) (  \phi \cdot ( (x'\#h') (y\# (\phi'\act g)) )  )\\
&=& \alpha(((x\#h)\#\phi) ( (x'\#h')(y\# (\phi'\act g)) ) \\
&=& \alpha(((x\#h)\#\phi) \circ \alpha((x'\#h')\#\phi') (y\#g),
\end{eqnarray*}
for all $x,x',y\in A,\,h,h',g\in H,\, \phi,\phi'\in H^*$,
therefore,  $\alpha$ is a homomorphism.
\end{proof}
\medskip

Let $\{ f_i\}$ be a basis of $H$ and $\{ \psi_i\}$ be the dual basis
of $H^*$, i.e., such that $\la f_i,\, \psi_j \ra =\delta_{ij}$ for
all $i,j$. Then we have identities
$$
\sum_i\, f_i \la h,\,\psi_i \ra = h,\qquad
\sum_i\, \la f_i ,\, \phi \ra \psi_i =\phi,
$$
for all $h\in H$ and $\phi\in H^*$, moreover the element
$\Sigma_i\, f_i\otimes \psi_i \in H\otimes H^*$ does not
depend on the choice of $\{ f_i\}$.

Let us define a linear map 
$\beta : \End(A\# H)_A \to (A\# H)\# H^*$ by
$$ 
\beta : T \mapsto \sum_i\, T(1\# {f_i}\2)(1\# S^{-1}({f_i}\1)) \# \psi_i.
$$

\begin{lemma}
\label{inverses}
The maps $\alpha$ and $\beta$ are inverses of each other.
\end{lemma}
\begin{proof}
We need to check that 
$$\beta\circ \alpha = \id_{(A\# H)\# H^*}
\quad \mbox{ and } \quad
\alpha \circ \beta  = \id_{\End(A\# H)_A}.
$$ 
For all $x\in A,\, h\in H$, and $\phi\in H^*$ we compute
\begin{eqnarray*}
\beta\circ \alpha((x\#h)\#\phi)
&=& \Sigma_i\,(x(h\1\cdot 1) \# h\2(\phi\act {f_i}\2)S^{-1}{f_i}\1) \# \psi_i\\
&=& \Sigma_i\,(x\# h \la \phi,\, {f_i}\3 \ra {f_i}\2 S^{-1}{f_i}\1) \# \psi_i\\
&=& \Sigma_i\,(x \# h\la \phi,\,1\2f_i\ra 1\1) \# \psi_i\\
&=& (x \# h (\phi\1\act 1))  \# \phi\2 \\
&=& (x \# h)  \# \eps_t(\phi\1) \phi\2  = (x \# h)\# \phi,
\end{eqnarray*}
where we used Lemma~\ref{a lemma} (iii) and the properties
of the element $\Sigma_i\,f_i\otimes \psi_i$.

Also, for every $T\in \End(A\# H)_A$ we have :
\begin{eqnarray*}
\alpha\circ\beta(T)(y\# g)
&=& \Sigma_i\, \alpha( T(1\# {f_i}\2)(1\# S^{-1}({f_i}\1)) \# \psi_i)(y\# g)\\
&=& \Sigma_i\, T(1\# {f_i}\2)(1\# S^{-1}({f_i}\1)) (y\# (\psi_i\act g) ) \\
&=& \Sigma_i\, T(1\# {f_i}\3)( (S^{-1}({f_i}\2) \cdot y) \# 
               S^{-1}({f_i}\1) g\1   ) \la \psi_i,\, g\2\ra \\
&=& T(1\# g\4)( (S^{-1}(g\3) \cdot y) \# S^{-1}(g\2) g\1) \\
&=& T(1\# g\3)( (S^{-1}(g\2) \cdot y) (\eps_s(g\1)\cdot 1) \# 1)\\
&=& T(1\# g\2) ((S^{-1}(g\1 1\2) \cdot y)(1\1 \cdot 1) \# 1) \\
&=& T(1\# g\2) ((S^{-1}(g\1) \cdot y)  \# 1) \\
&=& T( (g\2 S^{-1}(g\1) \cdot y) \# g\3) \\
&=& T((1\1\cdot y) \# 1\2g) =  T( y\# g),    
\end{eqnarray*}  
where we used that $T$ commutes with the right multiplication
by elements from $A$ and
identities from Lemma~\ref{a lemma}(i) and (iii).
\end{proof}

\begin{theorem}
\label{duality}
For any $H$-module algebra $A$ there is a canonical
isomorphism between the algebras
$(A\# H)\# H^*$ and $\End(A\# H)_A$.
\end{theorem}
\begin{proof}
Follows from Lemmas~\ref{alpha} and \ref{inverses}
\end{proof}

\begin{corollary}
$H\# H^* \cong \End(H)_{H_t}$, in particular, $H\# H^*$
is a semisimple algebra.
\end{corollary}
\begin{proof} 
We know that $H \cong H_t \# H$, where $H_t$ is the trivial 
$H$-module algebra, therefore applying Theorem~\ref{duality}
to $A=H_t$ we see that $H$ is a projective generating $H_t$-module
such that $\End(H)_{H_t} \cong H\# H^*$. Therefore, $H_t$ and
$H\# H^*$ are Morita equivalent. Since $H_t$ is always semisimple
(as a separable algebra), $H\# H^*$ is semisimple.
\end{proof}
\end{section}



\bibliographystyle{amsalpha}
  
\end{document}